\documentclass[a4paper,11pt]{article}
\usepackage{amsmath,amssymb,graphicx,epsf}

\newtheorem{thm}{\bfseries Theorem}
\newtheorem{lem}[thm]{\bfseries Lemma}        
\newtheorem{prop}[thm]{\bfseries Proposition} 
     
\newtheorem{defn}[thm]{\bfseries Definition}

\newenvironment{proof}{\medskip                    
\noindent{\scshape Proof:}}{\quad $\Box$\medskip}  

\def\NZQ{\mathbb}
\def\RR{{\NZQ R}}
\def\mA{{\mathcal A}}
\def\mB{{\mathcal B}}

\def\mM{{\mathcal M}}

\title{Every non-Euclidean oriented matroid admits a biquadratic final
 polynomial}
 
\date{October 24, 2005}

\author{Komei Fukuda%
\thanks{Research supported by the Swiss National Science Foundation Project  200021-105202,``Polytopes, Matroids and Polynomial Systems''.}\\ 
Swiss Federal Institue of Technology \\
Zurich and Lausanne, Switzerland\\
fukuda@ifor.math.ethz.ch
\and
Sonoko Moriyama\\
Graduate School of Information Science \\
and Technology, \\
The University of Tokyo, Japan\\
moriso@is.s.u-tokyo.ac.jp
\and
Hiroki Nakayama\\
Graduate School of Information Science \\ 
and Technology, \\
The University of Tokyo, Japan\\
nak-den@is.s.u-tokyo.ac.jp
}

\begin{document}
\maketitle


\begin{abstract}
 Richter-Gebert proved that every non-Euclidean uniform oriented matroid
 admits a biquadratic final polynomial. We extend this result to
 the non-uniform case. 
\end{abstract}


\section{Introduction}
In this paper, we identify a {\em chirotope\/}
 $\chi$ with an oriented matroid 
$\mM=(E,\chi)$, which we abbreviate by OM.  A standard reference for
the theory of oriented matroids is \cite{BLSWZ99}. 
The set $E=\{1,\ldots,n\}$ is called {\em ground
set\/} and $\chi:E^r\rightarrow\{+1,-1,0\}$ satisfies 
chirotope axioms,
where $r$ is a rank of an OM and $n$ is a number of elements
of the ground set. 

Let $X=(x_1,\ldots,x_n)\in\RR^{r\times n}$ be a 
configuration of $n$ points in $\RR^r$. Let $[i_1\cdots i_r]$ 
denote the determinant $\det(x_{i_1}\cdots x_{i_r})$.  By setting
$\chi_X(i_1,\ldots,i_r)={\rm sgn}[i_1\cdots i_r]$, the function $\chi_X$ 
satisfies the chirotope axioms.   A chirotope arising this way
is called {\em representable\/} or {\em realizable\/}.  
It is well known that not all chirotopes are realizable. 

In the sequel, we regard a bracket $[i_1\cdots i_r]$ as a
bracket variable. For any given ordered sequences of
indices $\tau=(\tau_1\cdots \tau_{r-2})$
and $\lambda=(\lambda_1\cdots\lambda_4)$, we call a bracket polynomial
\begin{eqnarray}
 [\tau\lambda_1\lambda_2][\tau\lambda_3\lambda_4] - 
  [\tau\lambda_1\lambda_3][\tau\lambda_2\lambda_4] +
  [\tau\lambda_1\lambda_4][\tau\lambda_2\lambda_3]
  \label{eq01}
\end{eqnarray}
a {\em 3-term Grassmann-Pl\"{u}cker polynomial}. If the chirotope is
realizable, the value of~(\ref{eq01}) is always $0$.
Now we introduce {\em biqudratic inequalities (equations)} and define
{\em biquadratic final polynomials}. 

\begin{defn}\normalfont
 \label{def01}
 Let $\chi$ be an OM of rank $r$, 
 let $\tau\in E^{r-2},\lambda\in E^4$ be
 index sequences, and let $A=(\tau\lambda_1\lambda_2)$,
 $B=(\tau\lambda_3\lambda_4)$, $C=(\tau\lambda_1\lambda_3)$,
 $D=(\tau\lambda_2\lambda_4)$, $E=(\tau\lambda_1\lambda_4)$ and
 $F=(\tau\lambda_2\lambda_3)$. Then
\begin{enumerate}
 \item  A pair $(\tau,\lambda)$ is called {\em $\chi$-normalized} if
       $\chi(A)\cdot\chi(B)\geq 0$, $\chi(C)\cdot\chi(D)\geq0$ and
       $\chi(E)\cdot\chi(F)\geq0$. 
  \item For a $\chi$-normalized pair $(\tau,\lambda)$, we call

       \begin{tabular}{cccc}
	$[A][B]<[C][D]$ & {\rm and} & $[E][F]<[C][D]$ & 
	{\em biquadratic inequalities}, \\
	$[A][B]=[C][D]$ & {\rm or} & $[E][F]=[C][D]$ & a {\em biquadratic equation}.
       \end{tabular} 
\end{enumerate} 

\end{defn}

We remark that for any pair $(\tau,\lambda)$, by permutating
$(\lambda_1,\lambda_2,\lambda_3,\lambda_4)$ appropriately,
$(\tau,\lambda)$ becomes $\chi$-normalized. We denote the set of
biquadratic inequalities and biquadratic equations by $\mA_\chi$ and
$\mB_\chi$, respectively. If $\chi$ is uniform, $\mB_\chi=\emptyset$.

\begin{defn}\normalfont
 \label{def02}
An OM $\chi$ admits a biquadratic final polynomial if
 there are a {\em non-empty} subset of $\mA\chi$ : $\{[A_i][B_i]<[C_i][D_i]
 \mid 1\leq i\leq k\}$ and a subset (maybe empty) of
 $\mB_\chi$ : $\{[A_j][B_j]=[C_j][D_j]\mid 1\leq j\leq l\}$ such that
 the following equality holds
 \begin{eqnarray*}
  \label{eq05}
   \prod_{i=1}^k[A_i][B_i]\cdot \prod_{j=1}^l[A_j][B_j] = 
   \prod_{i=1}^k[C_i][D_i]\cdot \prod_{j=1}^l[C_j][D_j]. 
 \end{eqnarray*}
\end{defn}
The following is a direct consequence of the definition above.

\begin{lem}
 \label{lm01}
 If $\chi$ admits a biquadratic final polynomial, $\chi$ is non-realizable.
\end{lem}

 Richter-Gebert~\cite{Ric93} proved that 
 every non-Euclidean uniform oriented matroid
 admits a biquadratic final polynomial. Our main theorem extends 
 this result to the non-uniform case.

\begin{thm} \label{maintheorem}
 Every non-Euclidean oriented matroid admits a biquadratic final
 polynomial. 
\end{thm}

\section{Oriented Matroid Programming}

Oriented matroid programming is formulated as a
combinatorial abstraction of linear programming~\cite{Bla77}. 
The simplex method in linear programming has a natural extension in
the setting of oriented matroids. Edmonds and Fukuda~\cite{Fuk82} showed
that there exist OMs allowing the simplex
method to generate a cycle of {\bf non-degenerate} pivots, which cannot%
 \footnote{Note that in linear programming, the simplex method can
generate a cycle of {\bf degenerate} pivots, known as cycling.}
occur in linear programming.
Consequently, one can show the non-realizability of
an OM by exhibiting a non-degenerate cycle of simplex pivots if 
exists.

Let $\chi$ be an OM of rank $r$ on an $(n+2)$ element set $E=\{1,\ldots,n,f,g\}$.  Here,
the last two elements $f$ and $g$ of $E$ are
distinguished.
The triple $(\chi,f,g)$ is called
an {\em oriented matroid program} (abbreviated by OMP). The element $g$ 
represents a hyperplane at infinity and $f$ represents an objective
function. 




\begin{defn}\normalfont
 \label{def03}
 Let $(\chi,f,g)$ be an OMP and $\mA$
 ($\mA^\infty$, respectively) be the affine (infinite) space with respect to $g$,
 i.e. the set of covectors with positive (zero) $g$-component.

\begin{enumerate}
 \item A set $B=(\lambda_1,\ldots,\lambda_{r-1})\in E-\{f,g\}$, such that
       $B\cup\{g\}$ is independent, is called an {\em affine basis}. The 
       unique vertex (i.e. a covector with minimal support, or 
       equivalently a cocircuit) $X$ with $X_B=0$ and $X_g=+$ is denoted by $v(B)$.
 \item $B_1\rightarrow B_2$ is called a {\em pivot operation} if
       $B_1,B_2$ are affine bases and $L=B_2-\{b\}=B_1-\{a\}$ where
       $a,b\in E-\{f,g\}$ and $a\neq b$. $L$ is called the {\em edge} of 
       $B_1\rightarrow B_2$.
 \item The {\em direction} of a pivot $L\cup\{a\}=B_1\rightarrow B_2=
       L\cup\{b\}$ where $L\cup\{a,b\}$ is assumed to be independent, is 
       the unique vertex $d=d(B_1\rightarrow B_2)\in\mA^\infty$ with
       $d_L=0$ and $d_a=v(B_2)_a$
 \item A pivot operation $L\cup\{a\}=B_1\rightarrow B_2=L\cup\{b\}$
       where $a\neq b$ is called\\
	{\em degenerate\/} if $v(B_1)=v(B_2)$, \\
	{\em horizontal\/}  if $L\cup\{f,g\}$ is dependent, \\
	{\em strictly increasing\/}  if $d(B_1\rightarrow B_2)_f>0$ and
		   $B_1\rightarrow B_2$ is not degenerate.
\end{enumerate}

\end{defn}

 We remark that neither degenerate nor horizontal pivot operation
 occurs  when an OM $\chi$ is uniform. 

\begin{defn}\normalfont
 \label{def04}
 A sequence of pivot operations $B_1\rightarrow B_2\rightarrow \cdots 
 \rightarrow B_k$ is called a {\em non-degenerate cycle} on $\chi$ if
 $B_1=B_k$ and all pivot operations are either degenerate, horizontal or 
 strictly increasing and at least one pivot is strictly increasing.
\end{defn}

Since no non-degenerate cycling occurs in linear programming, the
following proposition holds.
\begin{prop}
 If an OMP $(\chi,f,g)$ admits a non-degenerate cycle, then 
 the oriented matroid $\chi$ is
 non-realizable. 
\end{prop}

The following characterization of Euclidean OMs is fundamental.
\begin{prop}[\cite{Fuk82}]
 \label{lm02}
 An OMP $(\chi,f,g)$ on $E$ admits a non-degenerate cycle 
 for some choice of two distinguished elements $f$ and $g$ from $E$
 if and only if the oriented matroid $\chi$ is non-Euclidean. 
\end{prop}

\section{From Cycling to Biquadratic Final Polynomial} 

In the case of uniform OMs, Richter-Gebert~\cite{Ric93} 
gave a method to obtain a biquadratic final polynomial from a
non-degenerate cycle. 
Now we extend this method to the non-uniform case. In
the following proof, we translate each pivot operation to one
Grassmann-Pl\"{u}cker polynomial.

\begin{lem}
 \label{lm03}
 Let $(\chi,f,g)$ be an OMP and
 $L=\{\lambda_1,\ldots,\lambda_{r-2}\}\subset E-\{f,g\}$, $a,b\in E-\{f,g\}$ 
 such that $L\cup\{a\}=B_1\rightarrow B_2 =L\cup \{b\}$ is a pivot
 operation along edge $L$. Then 
\begin{itemize}
 \item if $B_1\rightarrow B_2$ is strictly increasing, \\
       $\chi(\lambda_1,\ldots,\lambda_{r-2},g,f)\cdot
       \chi(\lambda_1,\ldots,\lambda_{r-2},a,b)\cdot
       \chi(\lambda_1,\ldots,\lambda_{r-2},g,a)\cdot
       \chi(\lambda_1,\ldots,\lambda_{r-2},g,b)=+1$,
 \item if $B_1\rightarrow B_2$ is either degenerate or horizontal, \\
       $\chi(\lambda_1,\ldots,\lambda_{r-2},g,f)\cdot
       \chi(\lambda_1,\ldots,\lambda_{r-2},a,b)\cdot
       \chi(\lambda_1,\ldots,\lambda_{r-2},g,a)\cdot
       \chi(\lambda_1,\ldots,\lambda_{r-2},g,b)=0$.
\end{itemize}
\end{lem}

\begin{proof}
 For the first case, see~\cite{Ric93}. If the pivot operation is
 degenerate, which means two affine vertices 
 $v(B_1)$ and $v(B_2)$ are at same point,
 $\chi(\lambda_1,\ldots,\lambda_{r-2},a,b)=0$. Similary, if the pivot
 operation is horizontal, that is $L\cup\{f,g\}$ is dependent,
 $\chi(\lambda_1,\ldots,\lambda_{r-2},f,g)=0$ is satisfied. For both two 
 cases, the values become 0. 
\end{proof}

We are now ready to prove the main theorem.


\begin{proof} (of Theorem~\ref{maintheorem})
Let $\chi$ be a non-Euclidean OM on $E$.   By Proposition~\ref{lm02},
there exist $f$ and $g$ in $E$ such that the OMP $(\chi,f,g)$
admits a non-degenerate cycle, say, $B_1\rightarrow B_2\rightarrow \cdots \rightarrow B_k$ where $B_1=B_k$.
We shall construct a suitable biquadratic final polynomial. We define
 $L^i,a^i,b^i$ by the relations:
\begin{eqnarray}
 L^i\cup\{a^i\} = B_i\rightarrow B_{i+1}=L^i\cup\{b^i\} 
  \quad \mbox{\rm for all } 1\leq i\leq k.
  \label{eq06}
\end{eqnarray}
 In~(\ref{eq06}), we set
 $B_{k+1}=B_2$. $L^i=\{\lambda_1^i,\ldots,\lambda_{r-2}^i\}$ is the edge 
 of the pivot operation $B_i\rightarrow B_{i+1}$. We denote
 $\lambda^i=(\lambda_1^i,\ldots,\lambda_{r-2}^i)$.
 Consider the following sequence of Grassmann-Plu\"{u}cker polynomials: 
 \begin{eqnarray*}
  \label{eq07}
   GP^i= [\lambda^i,g,f][\lambda^i,a^i,b^i]-
   [\lambda^i,g,a^i][\lambda^i,f,b^i]+[\lambda^i,g,b^i][\lambda^i,f,a^i].
 \end{eqnarray*}

 Note that $GP^1=GP^k$. As in Definition~\ref{def01}, we set
 $A^i=(\lambda^i,g,f)$, 
 $B^i=(\lambda^i,a^i,b^i)$, $C^i=(\lambda^i,g,a^i)$,
 $D^i=(\lambda^i,f,b^i)$, $E^i=(\lambda^i,g,b^i)$ and
 $F^i=(\lambda^i,f,a^i)$. 
 Then, we have
 \[GP^i=A^i\cdot B^i-C^i\cdot D^i+E^i\cdot F^i.\]

Now we consider the signs of terms appearing in $GP^i$. If the pivot
 operation $B_i\rightarrow B_{i+1}$ is strictly increasing,
 $\chi(A^i)\cdot\chi(B^i)\cdot\chi(C^i)\cdot\chi(E^i)=+1$ is
 satisfied.  Using OM axioms,
the following 12 types of signs are possible:
 \[
 \begin{array}{cccccccccc}
  \label{tb01}
   A^i\cdot B^i & - & C^i &\cdot &D^i& + &E^i& \cdot & F^i & \\
  + & & + & & + & & + & & + & \mbox{{\rm type 1}} \\
  + & & + & & + & & + & & - & \mbox{{\rm type 2}} \\
  + & & + & & - & & + & & - & \mbox{{\rm type 3}} \\
  + & & - & & - & & - & & - & \mbox{{\rm type 4}} \\
  + & & - & & - & & - & & + & \mbox{{\rm type 5}} \\
  + & & - & & + & & - & & + & \mbox{{\rm type 6}} \\
  - & & + & & - & & - & & + & \mbox{{\rm type 7}} \\
  - & & + & & - & & - & & - & \mbox{{\rm type 8}} \\
  - & & + & & + & & - & & - & \mbox{{\rm type 9}} \\
  - & & - & & + & & + & & - & \mbox{{\rm type 10}} \\
  - & & - & & + & & + & & + & \mbox{{\rm type 11}} \\
  - & & - & & - & & + & & + & \mbox{{\rm type 12}} \\
 \end{array}
 \]
 After normalization, type 1, 4, 7 or 10 generates a biquadratic
 inequality $[E^i][F^i]<[C^i][D^i]$ and type 3, 6, 9 or 12 generates a
 biquadratic inequality $[A^i][B^i]<[C^i][D^i]$. 

 If the pivot operation $B_i\rightarrow B_{i+1}$ is degenerate or
 horizontal, $\chi(A^i)\cdot\chi(B^i)=0$ is satisfied. Using 
 OM axioms, the following 8 types of signs are possible.
 \[
 \begin{array}{cccccccccc}
  \label{tb02}
   A^i\cdot B^i & - & C^i &\cdot &D^i& + &E^i& \cdot & F^i & \\
  0 & & + & & + & & + & & + & \mbox{{\rm type 1'}} \\
  0 & & + & & + & & - & & - & \mbox{{\rm type 2'}} \\
  0 & & + & & - & & + & & - & \mbox{{\rm type 3'}} \\
  0 & & + & & - & & - & & + & \mbox{{\rm type 4'}} \\
  0 & & - & & + & & + & & - & \mbox{{\rm type 5'}} \\
  0 & & - & & + & & - & & + & \mbox{{\rm type 6'}} \\
  0 & & - & & - & & + & & + & \mbox{{\rm type 7'}} \\
  0 & & - & & - & & - & & - & \mbox{{\rm type 8'}} \\
 \end{array}
 \]
 Clearly, each one of the eight types implies a biquadratic equation $[E^i][F^i]=[C^i][D^i]$.

In both cases above, the fact
 $L^i\cup\{b^i\}=B_{i+1}=L^{i+1}\cup\{a^{i+1}\}$ implies the following
 relation: 
\begin{eqnarray*}
 \chi(D^i)\cdot\chi(E^i)\cdot\chi(C^{i+1})\cdot\chi(F^{i+1})=
  \chi(\lambda^i,f,b^i)\cdot\chi(\lambda^i,g,b^i)\cdot
  \chi(\lambda^{i+1},g,a^{i+1})\cdot
  \chi(\lambda^{i+1},f,a^{i+1}) \!&=&\! 1,
\end{eqnarray*}
 which restricts the types of possible successors $GP^{i+1}$ of a
 Grassmann-Pl\"{u}cker relation $GP^i$ of certain type. The transition 
 diagram is given in the following Figure~\ref{fig03}.

\begin{figure}[h]
 \begin{center}
  \includegraphics[height=45mm]{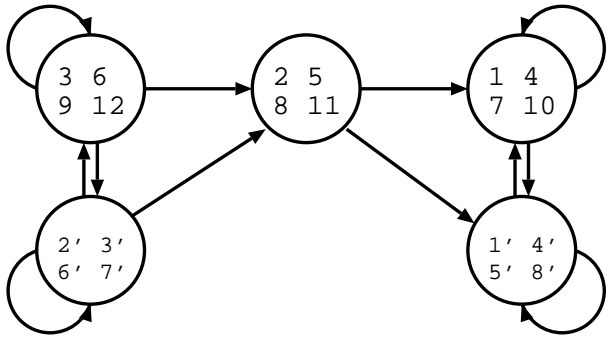}
  \caption{transition diagram among types}
  \label{fig03}  
 \end{center}
\end{figure}
 A Grassmann-Pl\"{u}cker relation of type $t$ can be succeeded by a 
 Grassmann-Pl\"{u}cker relation of type $s$ if and only if there is an
 arrow from the circle containing $t$ to the circle containing $s$. We
 have $GP^1=GP^k$ and $B_1\rightarrow B_2$ is strictly increasing, hence
 a sequence of transition is either
 \begin{itemize}
  \item contains only two states $(1,4,7,10)$ and $(1',4',5',8')$, and
	its initial state is $(1,4,7,10)$, or 
  \item contains only two states $(3,6,9,12)$ and $(2',3',6',7')$, and
	its initial state is $(3,6,9,12)$. 
 \end{itemize}
 In both cases, the resulting set of biquadratic inequalities and
 biquadratic equations yields a biquadratic final polynomial. 
\end{proof}





\end{document}